\documentclass[a4paper,12pt]{article}

\usepackage{amsmath}
\usepackage{amsthm}
\usepackage{amssymb}  
\usepackage{latexsym} 
\usepackage[all]{xy}
\usepackage{comment}
\usepackage{rotating}
\usepackage{slashbox}
\usepackage{multirow}
\usepackage{rotating}
\usepackage{enumerate}

\newcommand{\C}{{\mathbb C}}
\newcommand{\Fp}{{\mathbb{F}_p}}
\newcommand{\Fq}{{\mathbb{F}_q}}

\newfont{\wncyr}{wncyr10 at 12pt}
\newfont{\wncyrten}{wncyr10 at 10pt}

\newenvironment{Proof}{\par\noindent{\sc Proof:}}%
                      {\hspace*{\fill}\nobreak$\Box$\par\medskip}
                       {\hspace*{\fill}\nobreak$\Box$\par\medskip}
\newenvironment{myitemize}
{\begin{itemize}
\setlength{\itemsep}{1pt}
\setlength{\parskip}{0pt}
\setlength{\parsep}{0pt}}
{\end{itemize}}

\newtheorem{Proposition}{Proposition}[section]
\newtheorem{Theorem}[Proposition]{Theorem}
\newtheorem{Lemma}[Proposition]{Lemma}
\newtheorem{Corollary}[Proposition]{Corollary}

\theoremstyle{definition}

\newcounter{nootje}
\begin{document}
\normalsize
\title{Character Sums, Gaussian Hypergeometric Series, and a Family of Hyperelliptic Curves}
\author{Mohammad Sadek}
\date{}
\maketitle
\let\thefootnote\relax
\begin{abstract}{\footnotesize We study the character sums \[\phi_{(m,n)}(a,b)=\sum_{x\in\mathbb{F}_q}\phi\left(x(x^{m}+a)(x^{n}+b)\right),\textrm{ and, }
\psi_{(m,n)}(a,b)=\sum_{x\in\mathbb{F}_q}\phi\left((x^{m}+a)(x^{n}+b)\right)\] where $\phi$ is the quadratic character defined over $\Fq$. These sums are expressed in terms of Gaussian hypergeometric series over $\Fq$. Then we use these expressions to exhibit the number of $\Fq$-rational points on families of hyperelliptic curves and their Jacobian varieties.}
\end{abstract}
\section{Introduction}

Throughout this paper $p$ will denote an odd prime and $q=p^r$ where $r$ is a positive integer. We will write $\mathbb{F}_q$ for the finite field with $q$ elements. We let $\epsilon$ denote the trivial character and $\phi$ denote the quadratic character over $\mathbb{F}_q$. In other words, $\epsilon(x)=1$ for every $x\in\mathbb{F}_q^{\times}$, and $\epsilon(0)=0$; whereas $\phi(x)=1$ if $x$ is a square in $\mathbb{F}_q^{\times}$, $\phi(x)=-1$ if $x$ is not a square, and $\phi(0)=0$. If $a\in\mathbb{F}_q^{\times}$, the {\em Jacobsthal sum} $\phi_n(a)$ is defined by
\[\phi_n(a)=\sum_{x\in\mathbb{F}_q}\phi(x)\phi(x^n+a)\]and the modified Jacobsthal sum is defined by
\[\psi_n(a)=\sum_{x\in\mathbb{F}_q}\phi(x^n+a)\] where $n$ is a positive integer. Some of these sums are studied in literature, see for example \cite{Berndt} in which these sums are evaluated for small values of $n$. In \cite{Williams} the character sums
\[\sum_{x\in\mathbb{F}_p}\phi(ax^2+bx+c)\phi(Ax^2+Bx+C)\] are written in terms of simpler character sums. More precisely, the latter character sum is written in terms of a sum of the form $\displaystyle\sum_{x\in\mathbb{F}_p}\phi(f(x))$ where $f(x)$ is a cubic polynomial such that $y^2=f(x)$ is an elliptic curve. Consequently, some of these sums are evaluated explicitly which enables the author to find the number of rational points on elliptic curves defined via certain Weierstrass equations or quartic equations over prime fields.

Let $a,b\in\mathbb{F}_q$. In this note we are concerned with the following character sums
\begin{eqnarray*}
\phi_{(m,n)}(a,b)&=&\sum_{x\in\mathbb{F}_q}\phi(x)\phi(x^{m}+a)\phi(x^{n}+b),\textrm{ and}\\
\psi_{(m,n)}(a,b)&=&\sum_{x\in\mathbb{F}_q}\phi(x^{m}+a)\phi(x^{n}+b).
\end{eqnarray*}
These character sums can be thought of as generalizations of Jacobsthal sums and modified Jacobsthal sums respectively. We study the basic properties of these character sums. Then we focus on these sums when both $m$ and $n$ are powers of the prime $2$. In fact assuming that $q$ satisfies certain congruence relations we make use of elementary identities satisfied by Jacobi sums in order to express $\psi_{(m,n)}(a,b)$ and $\phi_{(m,n)}(a,b)$ in terms of specific Gaussian hypergeometric series. More precisely the character sums $\psi_{(m,n)}(a,b)$ and $\phi_{(m,n)}(a,b)$ are written as sums of values of Gaussian hypergeometric series ${_2 F_1}$ and ${_3F_2}$.

Since Greene initiated the study of hypergeometric series over finite fields in \cite{Greene}, many authors have written the number of rational points on algebraic curves over finite fields in terms of different Gaussian hypergeometric series, see for example \cite{BarmanKalita}, \cite{Barmanalgebraiccurves}, \cite{BarmanKalitaSaikia}, and \cite{Sadek}. Using identities relating the number of $\mathbb{F}_q$-rational points to hypergeometric series one can evaluate these series at specific values, see for example \cite{Ono}.

Given a hyperelliptic curve $C$ defined over $\Fq$ one can express the number of $\Fq$-rational points on $C$ using a sum of the quadratic character $\phi$ over elements in $\Fq$. In this note we are concerned with hyperelliptic curves whose number of $\Fq$-rational points is expressed using the sums $\psi_{(m,n)}(a,b)$ and $\phi_{(m,n)}(a,b)$. In particular, we study the number of $\Fq$-rational points on the hyperelliptic curves defined by the following affine equations
 \begin{eqnarray*}
 y^2&=&(x^m+a)(x^n+b),\;a,b\in\mathbb{F}_q^{\times},\;a\ne b,\\
 y^2&=&x(x^m+a)(x^n+b),\;a,b\in\mathbb{F}_q^{\times},\;a\ne b.
 \end{eqnarray*}
 Consequently, we establish a connection between the number of $\Fq$-rational points on these hyperelliptic curves and Gaussian hypergeometric series. One may go one step further and exploit the link between the number of $\Fq$-rational points on hyperelliptic curves and the number of $\Fq$-rational points on their Jacobians in order to produce explicit formulas for the latter number in terms of values of Gaussian hypergeometric series.

\section{Gaussian Hypergeometric Series}
In this section we introduce Gaussian hypergeometric series. Then we display some of the properties enjoyed by these series. We first recall that a {\em multiplicative character} $\chi$ over $\Fq$ is a homomorphism from $\Fq^{\times}\to \C^{\times}$, and we extend $\chi$ to $\Fq$ by setting $\chi(0) = 0$.

Let $J(A,B)$ denote the Jacobi sum \[J(A,B)=\sum_{x\in\Fq}A(x)B(1-x)\] where $A$ and $B$ are multiplicative characters over $\Fq$. Let $A_0,A_1,\ldots,A_n$ and $B_1,\ldots,B_n$ be multiplicative characters defined over $\Fq$. The {\em Gaussian hypergeometric series} is \[_{n+1}F_n\left(\begin{matrix}
                          A_0&A_1&\ldots&A_n \\
                          {} & B_1&\ldots&B_n
                        \end{matrix}\Big| x\right):=\frac{q}{q-1}\sum_{\chi}{{A_0\chi}\choose{\chi}}{{A_1\chi}\choose{B_1\chi}}\ldots{{A_n\chi}\choose{B_n\chi}}\chi(x)\]
where the sum is over all characters over $\Fq$, see Definition 3.10 in \cite{Greene}, and \[{A\choose B}:=\frac{B(-1)}{q}J(A,\overline{B})=\frac{B(-1)}{q}\sum_{x\in\mathbb{F}_q}A(x)\overline{B}(1-x)\]
and $\overline{B}$ is the character $1/B$, see Definition 2.4 in \cite{Greene}.

One can find the following properties of the symbol $\displaystyle{A\choose B}$  in \cite{Greene}.

\begin{Lemma}
\label{lem:charactersproperties}
Let $\epsilon$ and $\phi$ be the trivial and quadratic characters over $\mathbb{F}_q$ respectively. For any multiplicative characters $A$ and $B$ over $\mathbb{F}_q$, one has:
\begin{myitemize}
\item[a)] $\displaystyle A(1+x)=\delta(x)+\frac{q}{q-1}\sum_{\chi}{A\choose \chi}\chi(x)$ where $\delta(x)=1$ if $x=0$ and $\delta(x)=0$ if $x\ne 0$;
\item[b)] $\displaystyle \overline{A}(1-x)=\delta(x)+\frac{q}{q-1}\sum_{\chi}{A\chi\choose \chi}\chi(x)$ where $\delta(x)=1$ if $x=0$ and $\delta(x)=0$ if $x\ne 0$;
\item[c)] $\displaystyle {A\choose B}={A\choose A\overline{B}}$;
\item[d)] $\displaystyle {A\choose B}={B\overline{A}\choose B}B(-1)$;
\item[e)] $\displaystyle {A\choose \epsilon}={A\choose A}=-\frac{1}{q}+\frac{q-1}{q}\delta(A)$ where $\delta(A)=1$ if $A=\epsilon$ and $\delta(A)=0$ otherwise;
\item[f)] $\displaystyle {B^2\chi^2\choose \chi}={\phi B \chi\choose\chi}{B\chi\choose B^2\chi}{\phi\choose\phi B}^{-1}B\chi(4).$
\end{myitemize}
\end{Lemma}
\section{Evaluation of Some Character Sums}
In this section we link some character sums to certain hypergeometric series.
\begin{Lemma}
\label{lem1}
Let $f:\Fq\to\C$ be a map. One has
\[\sum_{x\in\mathbb{F}_q}\phi(x)f(x)=\sum_{x\in\mathbb{F}_q}f(x^2)-\sum_{x\in\mathbb{F}_q}f(x).\]
In particular one has $\psi_{2n}(a)=\psi_n(a)+\phi_n(a)$, and $\psi_{(2m,2n)}(a,b)=\psi_{(m,n)}(a,b)+\phi_{(m,n)}(a,b)$.
\end{Lemma}
\begin{Proof}
This holds since
\begin{eqnarray*}
\sum_{x\in\mathbb{F}_q}\phi(x)f(x)=\sum_{x\in\mathbb{F}_q}(1+\phi(x))f(x)-\sum_{x\in\mathbb{F}_q}f(x).
\end{eqnarray*}
The result follows as $1+\phi(x)$ is either $2$ if $x$ is a square in $\mathbb{F}_q^{\times}$; or $0$ otherwise.
\end{Proof}

\begin{Corollary}
\label{cor:jacobisum}
Let $a\in\mathbb{F}^{\times}_q$ and let $\chi, \psi$ be characters on $\mathbb{F}_q$. Then
\[\sum_{x\in\mathbb{F}_q}\psi(x^2)\chi(1+ax^2)=\psi\left(-a^{-1}\right)J(\psi,\chi)+\phi\psi\left(-a^{-1}\right)J(\phi\psi,\chi).\]
In particular if $\chi=\phi$ then
\[\sum_{x\in\mathbb{F}_q}\psi(x^2)\phi(1+ax^2)=q\phi(-1)\left[\psi\left(-a^{-1}\right){\psi\choose\phi\psi}+\phi\psi\left(-a^{-1}\right){\phi\psi\choose\psi}\right].\]
\end{Corollary}
\begin{Proof}
By setting $f(x)=\psi(x)\chi(1+ax)$, Lemma \ref{lem1} implies that
\begin{eqnarray*}
\sum_{x\in\Fq}\psi(x^2)\chi(1+ax^2)&=&\sum_{x\in\Fq}\psi(x)\chi(1+ax)+\sum_{x\in\Fq}\phi\psi(x)\chi(1+ax)
\end{eqnarray*}
Now, since the map $x\mapsto -a^{-1}x$ is bijective over $\mathbb{F}_q$, one has
\begin{eqnarray*}
\sum_{x\in\Fq}\psi(x^2)\chi(1+ax^2)&=&\psi\left(-a^{-1}\right)\sum_{x\in\Fq}\psi(x)\chi(1-x)+\phi\psi\left(-a^{-1}\right)\sum_{x\in\Fq}\phi\psi(x)\chi(1-x)\\
&=&\psi\left(-a^{-1}\right)J(\psi,\chi)+\phi\psi\left(-a^{-1}\right)J(\phi\psi,\chi),
\end{eqnarray*}
where the last equality holds from the definition of the Jacobi sum.
When $\chi=\phi$, one has
\begin{eqnarray*}
\sum_x\psi(x^2)\phi(1+ax^2)&=&\psi\left(-a^{-1}\right)J(\psi,\phi)+\phi\psi\left(-a^{-1}\right)J(\phi\psi,\phi)\\
&=&q\phi(-1)\left[\psi\left(-a^{-1}\right){\psi\choose\phi}+\phi\psi\left(-a^{-1}\right){\phi\psi\choose\phi}\right],
\end{eqnarray*}
where the last equality is implied by the definition of the symbol $\displaystyle{A\choose B}$.
Now one concludes using Lemma \ref{lem:charactersproperties} (c).
\end{Proof}
In fact one can generalize the identity above to reach more identities relating similar character sums to Jacobi sums. For example, if $q\equiv 1$ mod $4$ then there exists a character $\chi_4$ of order $4$ over $\Fq$. In this case if $\psi$ and $\chi$ are characters over $\Fq$ then Lemma \ref{lem1} implies that
\begin{eqnarray*}
\sum_{x\in\Fq}\psi(x^4)\chi(1+ax^4)&=&\sum_{x\in\Fq} \phi(x)\psi(x^2)\chi(1+ax^2)+\sum_{x\in\Fq}\psi(x^2)\chi(1+ax^2)\\
&=&\sum_{x\in\Fq}\chi_4\psi(x^2)\chi(1+ax^2)+\sum_{x\in\Fq}\psi(x^2)\chi(1+ax^2).
\end{eqnarray*}
Now using Corollary \ref{cor:jacobisum} one has
\begin{multline*}
\sum_{x\in\Fq}\psi(x^4)\chi(1+ax^4)= \chi_4\psi\left(-a^{-1}\right)J(\chi_4\psi,\chi)+\chi_4^3\psi\left(-a^{-1}\right)J(\chi_4^3\psi,\chi)\\+\psi\left(-a^{-1}\right)J(\psi,\chi)+\phi\psi\left(-a^{-1}\right)J(\phi\psi,\chi).
\end{multline*}
Therefore,
\begin{eqnarray*}
\sum_{x\in\Fp}\psi(x^4)\chi(1+ax^4)&=&\sum_{k=0}^3\chi_4^k\psi\left(-a^{-1}\right)J(\chi_4^k\psi,\chi).
\end{eqnarray*}
Now one can obtain the following result using a simple induction argument.
\begin{Theorem}
\label{thm:induction}
Let $m=2^t$ where $t$ is a positive integer such that $q\equiv 1$ mod $m$. Let $\chi, \psi$ be characters on $\mathbb{F}_q$, and $\chi_m$ a character of order $m$ over $\Fq$. Then one has
\begin{eqnarray*}
\sum_{x\in\Fq}\psi(x^m)\chi(1+ax^m)&=&\sum_{k=0}^{m-1}\chi_m^k\psi\left(-a^{-1}\right)J(\chi_m^k\psi,\chi).
\end{eqnarray*}
In particular, the following holds.
\begin{eqnarray*}
\sum_{x\in\Fq}\psi(x^m)\phi(1+ax^m)&=&q\phi(-1)\sum_{k=0}^{m-1}\chi_m^k\psi\left(-a^{-1}\right){\chi_m^k\psi\choose\chi_m^{k-\frac{m}{2}}\psi}.
\end{eqnarray*}
\end{Theorem}
\section{The Character Sums $\psi_{(m,n)}(a,b)$ and $\phi_{(m,n)}(a,b)$}
In this section we study the character sums $\psi_{(m,n)}(a,b)$ and $\phi_{(m,n)}(a,b)$. Moreover we express some of these sums using hypergeometric series.

Let $m$ and $n$ be nonnegative integers. Let $a$ and $b$ be elements in $\mathbb{F}_q^{\times}$. In this note we are studying the following character sums
\begin{eqnarray*}
\psi_{(m,n)}(a,b)&=&\sum_{x\in\Fq}\phi(x^m+a)\phi(x^n+b),\\
\phi_{(m,n)}(a,b)&=&\sum_{x\in\Fq}\phi(x)\phi(x^m+a)\phi(x^n+b).
\end{eqnarray*}
\subsection{Basic Properties}
 The following lemma collects some basic properties of these character sums.
\begin{Lemma}
\label{lem:sumsproperties}
Let $\displaystyle \chi_{(m,n)}\in\{\psi_{(m,n)},\phi_{(m,n)}\}$. The following statements hold.
\begin{enumerate}[i.]
\item $\displaystyle \chi_{(m,n)}(a,b)=\chi_{(n,m)}(b,a)$.
\item $\displaystyle\phi_{(n,0)}(a,b)=\phi(1+b)\phi_n(a)$, and, $\displaystyle\psi_{(n,0)}(a,b)=\phi(1+b)\psi_n(a)$.
\item $\displaystyle \psi_{(m,n)}(a^m,b)=\phi(a^{m+n})\psi_{(m,n)}\left(1,\frac{b}{a^n}\right)$.
\item $\displaystyle \phi_{(m,n)}(a^m,b)=\phi(a^{m+n+1})\phi_{(m,n)}\left(1,\frac{b}{a^n}\right)$.
\end{enumerate}
\end{Lemma}
\begin{Proof}
i) and ii) are obvious. For iii) one can verify the following equalities
\begin{eqnarray*}
\psi_{(m,n)}(a^m,b)&=&\sum_x \phi(x^m+a^m)\phi\left(x^n+b\right)=\phi(a^m)\sum_x \phi(x^m+1)\phi\left(a^nx^n+b\right)\\
&=& \phi\left(a^{m+n}\right)\sum_x\phi(x^m+1)\phi\left(x^n+\frac{b}{a^n}\right)
\end{eqnarray*}
where the second equality is obtained via replacing $x$ by $ax$. The proof of iv) is similar to that of iii).
\end{Proof}
\subsection{The Character Sum $\psi_{2,2}(a,b)$ and Elliptic Curves}
In \cite{Williams} the following identity was proved
\[\sum_{x\in\Fp}\phi\left(ax^2+bx+c\right)\phi\left(Ax^2+Bx+C\right)=\sum_{x\in\Fp}\left[\phi(x)\phi\left(Dx^2+\Delta x+d\right)\right]-\phi(aA)\]
where $a,b,c,A,B,C$ are integers; and $D,\Delta,d$ are defined by \[D=B^2-4AC,\;\Delta=4aC-2bB+4cA,\;d=b^2-4ac;\;\textrm{ and } \Delta^2-4Dd\not\equiv 0\textrm{ (mod } p).\]
The above equality yields the following formula for $\psi_{(2,2)}(a,b)$. \begin{eqnarray*}\psi_{(2,2)}(a,b)&=&\sum_{x\in\Fp}\phi(x^2+a)\phi(x^2+b)=\sum_{x\in\Fp}\phi\left(x(-4bx^2+4(a+b) x-4a)\right)-\phi(1)\\
&=&\phi(-b)\sum_{x\in\Fp}\phi\left(x(x-1)(x-ab^{-1})\right)-1.
\end{eqnarray*}
Given an elliptic curve $E$ defined over $\Fq$, the trace of the Frobenius of $E$ is $a_q(\lambda)=1+q-|E(\Fq)|$. In fact if $E$ is described by the Weierstrass equation $y^2=x(x-1)(x-\lambda)$ then $\displaystyle a_q(\lambda)=-\sum_{x\in\Fq}\phi\left(x(x-1)(x-\lambda)\right)$, see (6) in \cite{Ono}.
The following corollary is an immediate result of Theorem 1 of \cite{Ono}.
\begin{Corollary}
\label{cor22}
Let $a,b\in\mathbb{F}_p^{\times}$. Assuming that $a\ne b$, one has \[\psi_{(2,2)}(a,b)=-\phi(-b)a_p\left(ab^{-1}\right)-1=p\phi(b)\times {_2F_1}\left(\begin{matrix}
                          \phi&\phi \\
                          {} & \epsilon
                        \end{matrix}\Big| \frac{a}{b}\right)-1.\]
\end{Corollary}
Using the above corollary one can express the number of $\Fp$-rational points on the elliptic curve $E$ described by $y^2=(x^2+a)(x^2+b)$ as a hypergeometric series. This holds because $|E(\Fp)|=p+2+\psi_{(2,2)}(a,b)$.

In view of Theorem 2 of \cite{Ono}, the following result follows.
\begin{Corollary}
\label{cor22atvalues}
Let $a,b\in\mathbb{F}_p^{\times}$ be such that $a\ne b$. One has
\begin{align*}\psi_{(2,2)}(a,b)=\left\{\begin{array}{ll}
-1,  & \textrm{ if $p\equiv 3$ (mod $4$)} \\
2x\phi(b)(-1)^{(x+y+1)/2}-1, & \textrm{if $p\equiv 1$ (mod $4$), $p=x^2+y^2$, $x$ odd}, \end{array}\right.\end{align*}
where $b\in\{-a,2a,a/2\}$.
\end{Corollary}
In what follows we derive the same representation for $\psi_{(2,2)}(a,b)$ as a hypergeometric series if the base field is $\Fq$ instead of $\Fp$.
\begin{Proposition}
\label{proppsi22}
Let $a,b\in\mathbb{F}_q^{\times}$ be such that $a\ne b$. One has \[\psi_{(2,2)}(a,b):=\sum_{x\in\mathbb{F}_q}\phi(x^2+a)\phi(x^2+b)=q\phi(b)\times {_2F_1}\left(\begin{matrix}
                          \phi&\phi \\
                          {} & \epsilon
                        \end{matrix}\Big| \frac{a}{b}\right)-1.\]
\end{Proposition}
\begin{Proof}
Using Lemma \ref{lem:charactersproperties} a), one has
\begin{eqnarray*}
\psi_{(2,2)}(a,b)&=& \sum_{x\in\Fq}\phi(x^2+a)\phi(x^2+b)=\phi(ab)\sum_{x\in\Fq} \phi\left(1+\frac{x^2}{a}\right)\phi\left(1+\frac{x^2}{b}\right)\\
&=&\phi(ab)+\frac{q\phi(ab)}{q-1}\sum_{\chi}{\phi\choose\chi}\chi(b^{-1})\sum_{x\in\mathbb{F}^{\times}_q}\phi\left(1+\frac{x^2}{a}\right)\chi(x^2).
\end{eqnarray*}
The term $\phi(ab)$ appearing after the third equality corresponds to $x=0$.
Now Corollary \ref{cor:jacobisum} implies that
\begin{eqnarray*}
\psi_{(2,2)}(a,b)&=&\phi(ab)+\frac{q^2\phi(-ab)}{q-1}\sum_{\chi}{\phi\chi\choose\chi}\chi(-b^{-1})\left[\chi\left(-a\right){\chi\choose\phi\chi}+\phi\chi\left(-a\right){\phi\chi\choose\chi}\right]\\
&=& \phi(ab)+q\phi(-ab)\left[{_2F_1}\left(\begin{matrix}
                          \phi&\epsilon \\
                          {} & \phi
                        \end{matrix}\Big| \frac{a}{b}\right)+\phi(-a)\times{_2F_1}\left(\begin{matrix}
                          \phi&\phi \\
                          {} & \epsilon
                        \end{matrix}\Big| \frac{a}{b}\right)\right].
\end{eqnarray*}
Now one obtains
\begin{eqnarray*}
\psi_{(2,2)}(a,b)&=& \phi(ab)+q\phi(-ab)\left[{\phi\choose\phi}\phi\left(\frac{-a}{b}\right)-\frac{\phi\left(-1\right)}{q}\right]+q\phi(b)\times{_2F_1}\left(\begin{matrix}
                          \phi&\phi \\
                          {} & \epsilon
                        \end{matrix}\Big| \frac{a}{b}\right)\\
                        &=& -1+q\phi(b)\times{_2F_1}\left(\begin{matrix}
                          \phi&\phi \\
                          {} & \epsilon
                        \end{matrix}\Big| \frac{a}{b}\right),
\end{eqnarray*}
where the first equality follows from Corollary 3.16 ii) of \cite{Greene}, and the second equality holds since $\displaystyle{\phi\choose\phi}=\frac{-1}{q}$, see Lemma \ref{lem:charactersproperties} e), hence \begin{eqnarray*}q\phi(-ab)\left[{\phi\choose\phi}\phi\left(\frac{-a}{b}\right)-\frac{\phi\left(-1\right)}{q}\right]&=&
q\phi(-ab)\left[-\frac{1}{q}\phi\left(\frac{-a}{b}\right)-\frac{\phi\left(-1\right)}{q}\right]\\
&=&-\phi(1)-\phi(ab).\end{eqnarray*}
\end{Proof}

\subsection{The Character Sums $\phi_{(m,m)}(a,b)$ and $\psi_{(m,m)}(a,b)$}
Now we study the character sums $\psi_{(m,m)}(a,b)$ and $\phi_{(m,m)}(a,b)$ where $m$ is a positive integer, $q\equiv1$ (mod $2m$), and $a,b\in\mathbb{F}_q^{\times}$.
\begin{Lemma}
\label{lem:redefinepsi}
 Let $m$ be a positive integer such that $q\equiv 1$ (mod $2m$). Let $\chi_{2m}$ be a character of order $2m$ over $\Fq$. Then one has
 \begin{eqnarray*}
  \psi_{(m,m)}(a,b)&=& \sum_{x\in\Fq} \phi(x+a)\phi(x+b)\sum_{k=0}^{m-1}\chi_{2m}^{2k}(x),
 \end{eqnarray*}
 and
 \begin{eqnarray*}
 \phi_{(m,m)}(a,b)&=& \sum_{x\in\Fq} \phi(x+a)\phi(x+b)\sum_{k=0}^{m-1}\chi_{2m}^{2k+1}(x).
 \end{eqnarray*}
\end{Lemma}
\begin{Proof}
In order to show that the first equality holds, one uses the fact that $\chi_{2m}^2$ is a character of order $m$, therefore one has
\begin{align*}\sum_{k=0}^{m-1}\left(\chi_{2m}^{2}\right)^k(x)=\left\{\begin{array}{ll}
1  & \textrm{ if } x=0,\\
m & \textrm{ if  } x \textrm{ is an }m\textrm{-th power in }\mathbb{F}_q^{\times},\\
0 &\textrm{ otherwise.}
\end{array}\right.\end{align*}
Now recalling that $\displaystyle\psi_{(m,m)}(a,b)= \sum_{x\in\Fq}\phi(x^m+a)\phi(x^m+b)$, it follows that $\displaystyle\psi_{(m,m)}= \sum_{x\in\Fq} \phi(x+a)\phi(x+b)\sum_{k=0}^{m-1}\chi_{2m}^{2k}(x)$.

Now since $\displaystyle\psi_{(2m,2m)}=\sum_{x\in\Fq}\phi\left(x^{2m}+a\right)\phi\left(x^{2m}+b\right)$, and
\begin{align*}\sum_{k=0}^{2m-1}\chi_{2m}^k(x)=\left\{\begin{array}{ll}
1  & \textrm{ if } x=0,\\
2m & \textrm{ if  } x \textrm{ is a }2m\textrm{-th power in }\mathbb{F}_q^{\times},\\
0 &\textrm{ otherwise,}
\end{array}\right.\end{align*}
it follows that $\displaystyle \psi_{(2m,2m)}(a,b)= \sum_{x\in\Fq} \phi(x+a)\phi(x+b)\sum_{k=0}^{2m-1}\chi_{2m}^k(x)$.

The character sum  $\phi_{(m,m)}(a,b)$ is given by
\begin{eqnarray*}
\phi_{(m,m)}(a,b)=\psi_{(2m,2m)}(a,b)-\psi_{(m,m)}(a,b).
\end{eqnarray*}
Thus,
\begin{eqnarray*}
\phi_{(m,m)}(a,b)&=& \sum_{x\in\Fq} \phi(x+a)\phi(x+b)\sum_{k=0}^{2m-1}\chi_{2m}^k(x)- \sum_{x\in\Fq} \phi(x+a)\phi(x+b)\sum_{k=0}^{m-1}\chi_{2m}^{2k}(x)\\
&=&\sum_{x\in\Fq} \phi(x+a)\phi(x+b)\sum_{k=0,2\nmid k}^{2m-1}\chi_{2m}^k(x)\\
&=& \sum_{x\in\Fq} \phi(x+a)\phi(x+b)\sum_{k=0}^{m-1}\chi_{2m}^{2k+1}(x).
\end{eqnarray*}
\end{Proof}
Now one can express the character sums $\psi_{(m,m)}(a,b)$ and $\phi_{(m,m)}(a,b)$ in terms of Gaussian hypergeometric series of type ${_2F_1}$.
\begin{Theorem}
\label{corpsi22}
Let $m$ be a positive integer such that $q\equiv 1$ (mod $2m$). Let $\chi_{2m}$ be a character of order $2m$ over $\Fq$.
Then one has
\begin{eqnarray*}
\psi_{(m,m)}(a,b)&=&q\phi(-ab)\sum_{k=0}^{m-1}\chi_{2m}^{2k}(-a)\times{_2F_1}\left(\begin{matrix}
                          \phi&\chi_{2m}^{2k} \\
                          {} & \chi_{2m}^{2k+m}
                        \end{matrix}\Big| \frac{a}{b}\right),
\end{eqnarray*}
and
\begin{eqnarray*}\phi_{(m,m)}(a,b)&=&q\phi(-ab)\sum_{k=0}^{m-1}\chi_{2m}^{2k+1}(-a)\times{_2F_1}\left(\begin{matrix}
                          \phi&\chi_{2m}^{2k+1} \\
                          {} & \chi_{2m}^{2k+1+m}
                        \end{matrix}\Big| \frac{a}{b}\right).
                        \end{eqnarray*}
\end{Theorem}
\begin{Proof}
Lemma \ref{lem:redefinepsi} implies that
\begin{eqnarray*}
\psi_{(m,m)}(a,b)&=& \sum_{x\in\Fq} \phi(x+a)\phi(x+b)\sum_{k=0}^{m-1}\chi_{2m}^{2k}(x).
\end{eqnarray*}
One replaces $x$ with $-ax$ in the above equality in order to have
\begin{eqnarray*}
\psi_{(m,m)}(a,b)&=& \sum_{x\in\Fq} \phi(a-ax)\phi(b-ax)\sum_{k=0}^{m-1}\chi_{2m}^{2k}(-ax)\\
&=& \phi(ab)\sum_{k=0}^{m-1}\chi_{2m}^{2k}(-a)\sum_{x\in\Fq}\phi(1-x)\phi\left(1-\frac{a}{b}x\right)\chi_{2m}^{2k}(x)\\
&=& \phi(ab)\sum_{k=0}^{m-1}\chi_{2m}^{2k}(-a)\sum_{x\in\Fq}\chi_{2m}^{2k}(x) \left(\overline{\chi}_{2m}^{2k}\chi_{2m}^{2k+m}\right)(1-x)\phi\left(1-\frac{a}{b}x\right).
\end{eqnarray*}
Using \cite[Definition 3.5]{Greene}, the inner sum can be written in terms of a Gaussian hypergeometric series of type ${_2F_1}$
 \begin{eqnarray*}
 \psi_{(m,m)}(a,b)&=&\phi(ab)\sum_{k=0}^{m-1}\chi_{2m}^{2k}(-a)\frac{1}{\chi_{2m}^{4k+m}(-1)}\times q\times{_2F_1}\left(\begin{matrix}
                          \phi&\chi_{2m}^{2k} \\
                          {} & \chi_{2m}^{2k+m}
                        \end{matrix}\Big| \frac{a}{b}\right)\\
                        &=&q\phi(-ab)\sum_{k=0}^{m-1}\chi_{2m}^{2k}(-a)\times{_2F_1}\left(\begin{matrix}
                          \phi&\chi_{2m}^{2k} \\
                          {} & \chi_{2m}^{2k+m}
                        \end{matrix}\Big| \frac{a}{b}\right).
 \end{eqnarray*}
 The proof of the second equality is similar.
\end{Proof}

\section{$\psi_{(2,4)}$ and the Hypergeometric Series ${_3F_2}$}

\begin{Lemma}
\label{lem3f2}
Set $\lambda=a^2/b$. Assume that $\lambda\not\in\{0,-1\}$. Then one has
\[\left[1+\sum_{x\in\mathbb{F}_q}\phi\left(x(x+a)(x^2+b)\right)\right]^2=q+q^2\phi\left(\lambda+1\right)\times{_3F_2}\left(\begin{matrix}
                          \phi&\phi&\phi \\
                          {} & \epsilon&\epsilon
                        \end{matrix}\Big| \frac{\lambda}{\lambda+1}\right).\]
\end{Lemma}
\begin{Proof}
One has
\begin{eqnarray*}
\sum_{x\in\mathbb{F}_q}\phi\left(x(x+a)(x^2+b)\right)
&=& \sum_{x\in\mathbb{F}^{\times}_q}\phi(x^4)\phi\left(1+\frac{a}{x}\right)\phi\left(1+\frac{b}{x^2}\right)\\
&=&\sum_{x\in\mathbb{F}^{\times}_q}\phi\left(1+x\right)\phi\left(1+\frac{b}{a^2}x^2\right)\\
&=&\phi(b)\sum_{x\in\mathbb{F}^{\times}_q}\phi\left(1+x\right)\phi\left(\frac{a^2}{b}+x^2\right)\\
&=&\phi(-b)\sum_{x\in\mathbb{F}^{\times}_q}\phi\left(x-1\right)\phi\left(x^2+\frac{a^2}{b}\right).\\
\end{eqnarray*}
The trace of the Frobenius $a_q(\lambda)$ of the Clausen elliptic curve defined by the Weierstrass equation $y^2=(x-1)(x^2+\lambda)$ over $\Fq$ is given by $\displaystyle a_q(\lambda)=-\sum_{x\in\Fq}\phi(x-1)\phi(x^2+\lambda)$. The interested reader may consult \cite{OnoGuindy} for the arithmetic of rational points on Clausen elliptic curves over finite fields. Now it follows that
\begin{eqnarray*}
 \sum_{x\in\mathbb{F}_q}\phi\left(x(x+a)(x^2+b)\right)&=&-\phi(-b)\left[a_q(\lambda)+\phi(-\lambda)\right]\\
 &=&-\phi(-b)a_q(\lambda)-1.
 \end{eqnarray*}
 In other words, one has
\begin{eqnarray*}
a_q(\lambda)^2=\left[1+\sum_{x\in\mathbb{F}_q}\phi\left(x(x+a)(x^2+b)\right)\right]^2.
\end{eqnarray*}
 Now one uses Theorem 5 of \cite{Ono} which states that
\begin{eqnarray*}
a_q(\lambda)^2=q+q^2\phi\left(\frac{1}{\lambda+1}\right)\times{_3F_2}\left(\begin{matrix}
                          \phi&\phi&\phi \\
                          {} & \epsilon&\epsilon
                        \end{matrix}\Big| \frac{\lambda}{\lambda+1}\right).
\end{eqnarray*}
\end{Proof}
\begin{Theorem}
\label{thmpsi24}
Let $a,b\in\mathbb{F}_q^{\times}$. Set $\lambda=a^2/b$. We assume moreover that $\lambda\not\in\{0,-1\}$. Let $a_q(\lambda)$ be the trace of the Frobenius of the elliptic curve described by $y^2=(x-1)(x^2+\lambda)$. Then
\[\psi_{(2,4)}(a,b):=\sum_{x\in\Fq}\phi(x^2+a)\phi(x^4+b)=-\phi(-b)a_q(\lambda)-\phi(-a)a_q(1/\lambda)-1.\]
\end{Theorem}
\begin{Proof}
According to Lemma \ref{lem1}, one has
\begin{eqnarray*}
\sum_{x\in\mathbb{F}_q}\phi(x^2+a)\phi(x^4+b)=\sum_{x\in\mathbb{F}_q}\phi(x)\phi(x+a)\phi(x^2+b)+\sum_{x\in\mathbb{F}_q}\phi(x+a)\phi(x^2+b).
\end{eqnarray*}
We already saw during the course of the proof of Lemma \ref{lem3f2} that
\begin{eqnarray*}
\sum_{x\in\mathbb{F}_q}\phi(x)\phi(x+a)\phi(x^2+b)&=&-\phi(-b)a_q(\lambda)-1.
\end{eqnarray*}
Furthermore, by replacing $x$ by $-ax$, one knows that \[ \sum_{x\in\mathbb{F}_q}\phi(x+a)\phi(x^2+b)=\phi(-a)\sum_{x\in\mathbb{F}_q}\phi(x-1)\phi\left(x^2+\frac{b}{a^2}\right)=-\phi(-a)a_q(1/\lambda).\]
It follows that
\begin{eqnarray*}
\sum_{x\in\mathbb{F}_q}\phi(x^2+a)\phi(x^4+b)&=&
 -\phi(-b)a_q(\lambda)-\phi(-a)a_q(1/\lambda)-1.
\end{eqnarray*}
\end{Proof}
Now one can express the character sum $\psi_{(2,4)}(a,b)$ in terms of the Gaussian hypergeometric series $\displaystyle {_3F_2}\left(\begin{matrix}
                          \phi&\phi&\phi \\
                          {} & \epsilon&\epsilon
                        \end{matrix}\Big| \frac{\lambda}{\lambda+1}\right)$, where $\lambda=a^2/b$, up to a sign ambiguity, see Lemma \ref{lem3f2}.
\begin{Corollary}
\label{corpsi24}
Let $a,b\in\mathbb{F}_q^{\times}$. Set $\lambda=a^2/b$. Assume moreover that $\lambda\not\in\{0,-1\}$. Let $q\equiv 1$ (mod $4$) and $\chi_4$ a character of order $4$ defined over $\Fq$. Then
\[\psi_{(2,4)}(a,b)=q\phi(-b)\times{_2F_1}\left(\begin{matrix}
                          \chi_4&\chi_4^3 \\
                          {} & \epsilon
                        \end{matrix}\Big| -\lambda\right)+q\phi(-a)\times{_2F_1}\left(\begin{matrix}
                          \chi_4&\chi_4^3 \\
                          {} & \epsilon
                        \end{matrix}\Big| \frac{-1}{\lambda}\right)-1.\]
\end{Corollary}
\begin{Proof}
Proposition 1.2 in \cite{OnoGuindy} indicates that
\[a_q(\lambda)=-q\times{_2F_1}\left(\begin{matrix}
                          \chi_4&\chi_4^3 \\
                          {} & \epsilon
                        \end{matrix}\Big| -\lambda\right).\]
                        Now the result follows using Theorem \ref{thmpsi24}.
\end{Proof}

\section{Hyperelliptic curves and their Jacobians}

We may express the number of $\mathbb{F}_q$-rational points on a certain hyperelliptic curve in terms of a character sum involving the quadratic character $\phi$. This holds because a hyperelliptic curve $C$ is the union of two affine pieces
\[y^2=f(x),\textrm{ and } z^2=  x^{2g+2}f(1/x)\] where $f(x)\in\Fq[x]$. If $C$ has an $\Fq$-rational Weierstrass point then one may assume that $\deg f(x)$ is odd. In fact the genus of the curve $C$ is $g$ if $\deg f(x)$ is either $2g+1$ or $2g+2$.
In this section the polynomial $f(x)$ is monic. Thus one can think of $C$ as the
smooth projective model of the affine curve $y^2=f(x)$, which is obtained by adding one or two points at infinity according to $\deg f(x)$ being odd or even respectively.
In particular one has
\begin{eqnarray}
\label{eq1}
|C(\Fq)|&=& r+\#\{(x,y)\in\mathbb{F}^2_q: y^2=f(x)\}\nonumber\\
&=&r+\sum_{x\in\Fq}\left[1+\phi(f(x))\right]\nonumber\\
&=&r+q+\sum_{x\in\Fq}\phi(f(x))
\end{eqnarray}
where $r=1$ if $\deg f(x)$ is odd; and $r=2$ if $\deg f(x)$ is even, see \cite[Ch V, \S 1]{sil1} for this equation when $f(x)$ is a polynomial of degree $3$. In other words $r$ is the number of rational points at infinity.

One then can use Hasse-Weil bound on the size of the set of rational points over $\Fq$ in order to find a bound on the absolute value of the corresponding character sum. One recalls that if $C$ is an irreducible nonsingular algebraic curve over $\Fq$ of genus $g\ge1$, then Hasse-Weil Theorem gives that
\[\Big||C(\Fq)|-q-1\Big|\le 2g \sqrt{q}.\]

Furthermore one can use the following lemma to evaluate the number of rational points on the Jacobian varieties of several hyperelliptic curves over a prime field. For the definition of the Jacobian of an abelian variety, the reader may consult \cite[Definition 7.4.40]{Liu}.

\begin{Lemma}
\label{lemjacobian}
Let $C$ be a smooth projective curve of genus $g> 1$. Let $J_C$ be the Jacobian of $C$. If $N_k=|C(\mathbb{F}_{p^k})|$ then
\begin{align*}|J_C(\Fp)|=\left\{\begin{array}{ll}
\frac{N_1^2+N_2}2 -p & \textrm{ if } g=2 \\
\frac{N_1^3}6+\frac{N_1N_2}2+\frac{N_3}3-pN_1 & \textrm{ if }g=3\\
1+p^g&\textrm{ if } N_k=1+p^k,\;k=1,2,\ldots,g. \end{array}\right.\end{align*}
\end{Lemma}

In \cite{Berndt} the authors evaluated the Jacobsthal sums $\psi_n(a)$ and $\phi_n(a)$ for small values of $n$. For the convenience of the reader we are going to write the number of $\Fp$-rational points on the hyperelliptic curves described by the affine equations $y^2=x^n+a$ or $y^2=x(x^n+a)$ for some values $n\le 8$ using the results in \cite{Berndt} and the equality given in (\ref{eq1}).

\begin{Proposition}
\label{prop1}
Let $p$ be a prime. Let $C_n$ be a hyperelliptic curve defined by the affine equation $y^2=x^n+a$, and $C_n'$ defined by $y^2=x(x^n+a)$ where $p\nmid a$. Then the following statements hold.
\begin{enumerate}[i.]
\item  Assume $p \equiv 1$ (mod $6$) where $p=a_3^2+3b_3^2$ and $a_3\equiv -1$ (mod $3$). Let $\eta_3=\pm1$. One has
\begin{align*}|C_3'(\Fp)|=\left\{\begin{array}{ll}
1+p+2\phi(a)a_3,  & \textrm{if $a$ is a cubic (mod }p) \\
1+p-\phi(a)(a_3+3\eta_3|b_3|), & \textrm{otherwise}, \end{array}\right.\end{align*}

\begin{align*}|C_6(\Fp)|=\left\{\begin{array}{ll}
1+p+2\left[1+\phi(a)\right]a_3,  & \textrm{if $a$ is a cubic (mod }p) \\
1+p-\phi(a)(a_3+3\eta_3|b_3|)-(a_3-3\eta_3|b_3|), & \textrm{otherwise}, \end{array}\right.\end{align*}

\begin{align*}|C_3'(\Fp)|=\left\{\begin{array}{ll}
1+p+2a_3,  & \textrm{if $a$ is a cubic (mod }p) \\
1+p-(a_3-3\eta_3|b_3|), & \textrm{otherwise}. \end{array}\right.\end{align*}
\item Assume $p\equiv 1$ (mod $4$) where $p=c_4^2+d_4^2$ and $c_4\equiv -1$ (mod $4$). One has
\begin{align*}|C_4(\Fp)|=|C_2'(\Fp)|=\left\{\begin{array}{ll}
1+p+2c_4,  & \textrm{if $a$ is a quartic (mod }p) \\
1+p-2c_4, & \textrm{if $a$ is a quadratic but not a quartic}\\
1+p\pm2|d_4| &\textrm{otherwise,}
 \end{array}\right.\end{align*}
 \item Assume that $p\equiv 1$ (mod $8$) where $p=a_8^2+2b_8^2$ and $a_8\equiv -1$ (mod $4$). One has
 \begin{align*}|C_4'(\Fp)|=\left\{\begin{array}{ll}
1+p+4(-1)^{(p-1)/8}a_8,  & \textrm{if $a$ is an octic (mod }p) \\
1+p-4(-1)^{(p-1)/8}a_8, & \textrm{if $a$ is a quartic but not an octic (mod } p)\\
1+p &\textrm{if $a$ is a quadratic but not a quartic (mod }p)\\
1+p\pm 4|b_8| &\textrm{otherwise,}
 \end{array}\right.\end{align*}
 \begin{align*}|C_8(\Fp)|=\left\{\begin{array}{ll}
1+p+2c_4+4(-1)^{(p-1)/8}a_8,  & \textrm{if $a$ is an octic (mod }p) \\
1+p+2c_4-4(-1)^{(p-1)/8}a_8, & \textrm{if $a$ is a quartic but not an octic (mod } p)\\
1+p-2c_4 &\textrm{if $a$ is a quadratic but not a quartic (mod }p)\\
1+p\pm 2|d_4|\pm 4|b_8| &\textrm{otherwise.}
 \end{array}\right.\end{align*}
\end{enumerate}
\end{Proposition}

Using the results in \S 4 of \cite{Berndt} one can find $|C_6'(\Fp)|$ and $|C_{12}(\Fp)|$ for $p\equiv 1$ (mod $12$); $|C_{12}'(\Fp)|$ and $|C_{24}(\Fp)|$ for $p\equiv 1$ (mod $24$). In fact one can obtain the number of $\mathbb{F}_{p^2}$-rational points on the same curves using similar results over $\mathbb{F}_{p^2}$ in \cite{BerndtIllinois}.

In \cite{NIITSUMA} explicit formulas for the number of $\Fp$-rational points on hyperelliptic curves defined by $y^2=x^6+a$, $y^2=x(x^6+a)$, and $y^2=x^{12}+a$ are given.

Since we expressed some character sums $\psi_{(m,n)}(a,b)$ and $\phi_{(m,n)}(a,b)$ in terms of Gaussian hypergeometric series, and knowing that these character sums are relevant to the number of rational points on hyperelliptic curves, see (\ref{eq1}), one will represent the number of rational points using Gaussian hypergeometric series.

\begin{Theorem}
 Let $m$ be a positive integer such that $q\equiv 1$ (mod $2m$). Let $\chi$ be a character of order $2m$ over $\Fq$. Let $C_m$ and $C'_m$ be hyperelliptic curves defined by the affine equations $y^2=(x^m+a)(x^m+b)$ and $y^2=x(x^{m}+a)(x^{m}+b)$ over $\Fq$, respectively.
  One has
 \begin{eqnarray*}|C_m(\Fq)|&=&2+q+\psi_{(m,m)}(a,b)\\&=&2+q+q\phi(-ab)\sum_{k=0}^{m-1}\chi^{2k}(-a)\times{_2F_1}\left(\begin{matrix}
                          \phi&\chi^{2k} \\
                          {} & \chi^{2k+m}
                        \end{matrix}\Big| \frac{a}{b}\right),
 \end{eqnarray*}
 and
 \begin{eqnarray*}|C'_m(\Fq)|&=&1+q+\phi_{(m,m)}(a,b)\\&=&1+q+q\phi(-ab)\sum_{k=0}^{m-1}\chi^{2k+1}(-a)\times{_2F_1}\left(\begin{matrix}
                          \phi&\chi^{2k+1} \\
                          {} & \chi^{2k+1+m}
                        \end{matrix}\Big| \frac{a}{b}\right).\end{eqnarray*}
\end{Theorem}
\begin{Proof}
This is Theorem \ref{corpsi22} and (\ref{eq1}).
\end{Proof}

\begin{Theorem}
\label{thmcurvepsi24}
Let $a,b\in\mathbb{F}_q^{\times}$. Set $\lambda=a^2/b$. We assume moreover that $\lambda\not\in\{0,-1\}$. Let $a_q(\lambda)$ be the trace of the Frobenius of the elliptic curve described by $y^2=(x-1)(x^2+\lambda)$. Let $C$ be the hyperelliptic curve defined by the affine equation $y^2=(x^2+a)(x^4+b)$. Then
\begin{eqnarray*}
|C(\Fq)|&=&2+q+\psi_{(2,4)}(a,b)=1+q-\phi(-b)a_q(\lambda)-\phi(-a)a_q(1/\lambda).
\end{eqnarray*}
In particular if $q\equiv 1$ (mod $4$) and $\chi$ a character of order $4$ over $\Fq$ then
\[|C(\Fq)|=1+q+q\phi(-b)\times{_2F_1}\left(\begin{matrix}
                          \chi&\chi^3 \\
                          {} & \epsilon
                        \end{matrix}\Big| -\lambda\right)+q\phi(-a)\times{_2F_1}\left(\begin{matrix}
                          \chi&\chi^3 \\
                          {} & \epsilon
                        \end{matrix}\Big| \frac{-1}{\lambda}\right).\]
                        Moreover one has \[\left|\phi(-b)\times{_2F_1}\left(\begin{matrix}
                          \chi&\chi^3 \\
                          {} & \epsilon
                        \end{matrix}\Big| -\lambda\right)+\phi(-a)\times{_2F_1}\left(\begin{matrix}
                          \chi&\chi^3 \\
                          {} & \epsilon
                        \end{matrix}\Big| \frac{-1}{\lambda}\right)\right|\le\frac{4}{\sqrt q}.\]
 \end{Theorem}
\begin{Proof}
This is Theorem \ref{thmpsi24}, Corollary \ref{corpsi24} and (\ref{eq1}). The estimate is obtained using Hasse-Weil bound.
\end{Proof}
Since the hyperelliptic curve $C$ in Theorem \ref{thmcurvepsi24} is of genus 2 one may use Lemma \ref{lemjacobian} to express the number of $\Fp$-rational points on the Jacobian $J$ of $C$ in terms of hypergeometric series.

\bibliographystyle{plain}
\footnotesize
\bibliography{Edwards}
Department of Mathematics and Actuarial Science\\ American University in Cairo\\ mmsadek@aucegypt.edu
\end{document}